\documentclass{amsart}

\usepackage[ngerman,american]{babel}
\usepackage{graphicx,color}
\usepackage{a4wide}
\usepackage{amssymb}
\usepackage[latin1]{inputenc}
\usepackage[hyperref]{hyperref}
\usepackage{psfig}

\renewcommand{\phi}{\varphi}
\newcommand{\C}{{\mathbb{C}}}
\newcommand{\R}{{\mathbb{R}}}
\newcommand{\Z}{{\mathbb{Z}}}
\newcommand{\N}{{\mathbb{N}}}

\renewcommand{\epsilon}{\varepsilon}
\renewcommand{\theta}{\vartheta}
\renewcommand{\S}{{\mathbb{S}}}
\newcommand{\Disc}[1]{{\mathbb{D}^{#1}}}
\newcommand{\OpenDisc}[1]{{{\mathbb{D}}^{#1}_{<1}}}
\newcommand{\RP}[1]{{\R\mathbb{P}^{#1}}}
\newcommand{\CP}[1]{{\C\mathbb{P}^{#1}}}

\newcommand{\g}{{\mathfrak{g}}}

\newcommand{\lie}[1]{{\mathcal{L}_{#1}}}
\newcommand{\pairing}[2]{{\langle{#1}|{#2}\rangle}}
\newcommand{\norm}[1]{{\lVert #1\rVert}}
\newcommand{\abs}[1]{{\lvert #1\rvert}}
\newcommand{\intersection}[2]{{\iota(#1,#2)}}
\newcommand{\regorb}[1]{{{#1}_{(\mathrm{reg})}}}
\newcommand{\princorb}[1]{{{#1}_{(\mathrm{princ})}}}
\newcommand{\singorb}[1]{{{#1}_{(\mathrm{sing})}}}
\newcommand{\Dehn}{\tau}
\newcommand{\POS}{{\mathbf{q}}}
\newcommand{\MOM}{{\mathbf{p}}}
\newcommand{\lcan}{{\lambda_{\mathrm{can}}}}
\newcommand{\Etriv}{{E_{\mathrm{triv}}}}
\newcommand{\Etwist}{{E_{\mathrm{twist}}}}

\DeclareMathOperator{\Ad}{Ad}
\DeclareMathOperator{\ad}{ad}

\DeclareMathOperator{\id}{id}
\DeclareMathOperator{\orb}{Orb}

\DeclareMathOperator{\rot}{rot}
\DeclareMathOperator{\SO}{SO}
\DeclareMathOperator{\orthgroup}{O}
\newcommand{\so}[1]{{\mathfrak{so}({#1})}}
\DeclareMathOperator{\stab}{Stab}
\DeclareMathOperator{\SU}{SU}

\theoremstyle{plain}
\newtheorem{theorem}{Theorem}
\newtheorem{lemma}[theorem]{Lemma}
\newtheorem{coro}[theorem]{Corollary}

\theoremstyle{remark}
\newtheorem{remark}{Remark}
\newtheorem{example}{Example}
\newtheorem*{brieskornbeispiel}{Example~\ref{brieskorn beispiel}}
\newtheorem*{sphaerenbeispiel}{Example~\ref{sphaeren beispiel}}

\theoremstyle{definition}
\newtheorem*{defi}{Definition}

\begin{document}
\bibliographystyle{amsalpha}

\title[$5$-dimensional contact $\SO(3)$-manifolds]{$5$-dimensional
  contact $\SO(3)$-manifolds and Dehn twists}

\author{Klaus Niederkrüger}

\address{Mathematisches Institut, Universität zu Köln\\Weyertal
  86-90\\50.931 Köln\\Federal Republic of Germany}

\email{kniederk@math.uni-koeln.de}

\begin{abstract}
  In this paper the $5$-dimensional contact $\SO(3)$-manifolds are
  classified up to equivariant contactomorphisms.  The construction of
  such manifolds with singular orbits requires the use of generalized
  Dehn twists.
  
  We show as an application that all simply connected $5$-manifolds
  with singular orbits are realized by a Brieskorn manifold with
  exponents $(k,2,2,2)$.  The standard contact structure on such a
  manifold gives right-handed Dehn twists, and a second contact
  structure defined in the article gives left-handed twists.
\end{abstract}

\maketitle

A $5$-dimensional contact $\SO(3)$-manifold $M$ can be decomposed into
the set of singular orbits $\singorb{M}$ and the set of regular orbits
$\regorb{M}$. Both parts can be described relatively easily: The
singular orbits are the disjoint union of copies of
$\S^1\times\RP{2}$, $\S^1\times\S^2$ or
$\S^1\widetilde\times\S^2:=\R\times\S^2/\!\sim$, where
$(t,p)\sim(t+1,-p)$.  The set of regular orbits contains a canonical
submanifold $R$ of dimension $3$ (the so-called cross-section), and
one has that $\regorb{M}\cong\SO(3)\times_{\S^1}R$.

For gluing the singular orbits onto the regular ones, there is an
integer invariant that classifies all possibilities. This integer
corresponds to the number of Dehn twists.

\subsection*{Acknowledgments}
This article contains most of the results of my Ph.D. thesis. I would
like to thank the \emph{Universiteit Leiden} and the \emph{Universität
  zu Köln} for funding my work.  I am indebted to my supervisor
Hansjörg Geiges for his support and patience.  Furthermore I'm
grateful to the following people for many useful discussions: Peter
Heinzner, Federica Pasquotto, Juan Souto, and Kai Zehmisch, but most
of all to Otto van Koert.

\setcounter{section}{-1}
\section{Notation}\label{notation}

This section only fixes some notations about Lie groups and
$G$-manifolds.  In the article, $G$ denotes always a compact,
connected Lie group, $\g$ is its Lie algebra and $\g^*$ is the
corresponding coalgebra.  The only $G$-operation considered on $\g^*$
will be the coadjoint action.  For the stabilizer of an element
$\nu\in\g^*$, we write $G_\nu$.

A $G$-equivariant map $\Phi$ between $G$-manifolds $M$ and $N$
consists of a smooth map $\Phi_M:\,M\to N$ such that $\Phi_M(gp)=
g\,\Phi_M(p)$.  As a short-hand, we will write $G$-diffeomorphism
instead of $G$-equivariant diffeomorphism, $G$-contactomorphism
instead of $G$-equivariant contactomorphism etc.

Let $N$ be a submanifold of a $G$-manifold $M$.  The \textbf{flow-out
  of $N$} is defined as the set $G\cdot N$.

For a $G$-manifold $M$, we denote the set of principal orbits by
$\princorb{M}$, the set of singular orbits by $\singorb{M}$ and the
set of regular (i.e.~non-singular) orbits by $\regorb{M}$.  The
conjugation class of a closed subgroup $H\le G$ is written $(H)$, and
$M_{(H)}$ is the set of points $p\in M$ whose stabilizer $\stab(p)$
lies in the class $(H)$.  The normalizer $N(H)$ of $H$ is the subgroup
$\{g\in G|gHg^{-1}=H\}$.

For every element $X\in\g$, the \textbf{infinitesimal generator} of the
action is the vector field $X_M(p) :=
\left.\frac{d}{dt}\right|_{t=0}\exp(tX)p$.

\section{Preliminaries}

At any point $p$ of a $G$-manifold, there exists a so-called
\textbf{slice} $S_p$.  This is a submanifold that is transverse to the
orbit $\orb(p)$, invariant under the action of $\stab(p)$, and
satisfies the condition that whenever $g\cdot q\in S_p$ (with $g\in G$
and $q\in S_p$), then $g\in\stab(p)$. In particular for the coadjoint
action on $\g^*$, there exists a unique maximal slice at any
$\nu\in\g^*$, which will be denoted by $S_\nu^*$
(see~\cite{Duistermaat}).

\begin{example}\label{SO(3)-struktur}
  Consider the $\SO(3)$-structure of $\so{3}^*$ given by the coadjoint
  action.  The principal orbits are $2$-spheres lying concentrically
  around $0$, and $\{0\}$ is the only singular orbit in $\so{3}^*$.
  The maximal slice of an element $\nu\in\so{3}^*$ ($\nu\ne 0$) is
  $\R^+\cdot\nu$ and the maximal slice at $0$ is the whole of
  $\so{3}^*$.
\end{example}

\begin{defi}
  A \textbf{contact $G$-manifold} $(M,\alpha)$ is a $G$-manifold with
  an invariant contact form $\alpha$.\footnote{If $\xi=\ker\alpha$ is
    a $G$-invariant contact structure on $M$, then one can average
    $\alpha$ over the $G$-action to obtain an invariant contact form.}
\end{defi}

In the rest of the article we will assume that all contactomorphisms
\textbf{preserve the coorientation} of the contact structure, i.e.\ 
for a contactomorphism $\Phi:\,(M,\alpha)\to (M^\prime,\alpha^\prime)$
with $\Phi^*\alpha^\prime = f\,\alpha$, the function $f$ has to be
positive.

\begin{defi}
  The \textbf{moment map} $\mu:\,M\to\g^*$ of a contact $G$-manifold
  $(M,\alpha)$ is given by
  $$
  \pairing{\mu(p)}{X} := \alpha_p(X_M)\;.
  $$
\end{defi}

\begin{defi}
  For a contact $G$-manifold $(M,\alpha)$ with moment map
  $\mu:\,M\to\g^*$, the \textbf{cross-section $R$ at a point
    $\nu\in\mu(M)$} is defined as
  $$
  R := \mu^{-1}(S_\nu^*)\;.
  $$
\end{defi}

One can find a symplectic version of the following theorem in
\cite{Lerman_Meinrenken}, the contact version has been described in
\cite{Willett}.

\begin{theorem}[cross-section theorem]\label{cross-section}
  Let $(M,\alpha)$ be a contact $G$-manifold with moment map
  $\mu_M:\,M\to\g^*$.  Let $\nu\in\g^*$ be an element in the image of
  the moment map, and let $S_\nu^*\subseteq\g^*$ be the unique maximal
  slice at $\nu$.
  
  Then:
  \begin{enumerate}
  \item The \textbf{cross-section} $R:=\mu_M^{-1}(S_\nu^*)$ is a contact
    $G_\nu$-submanifold of $M$, where $G_\nu:=\stab(\nu)$.
  \item The $G$-action induces a $G$-diffeomorphism between the
    flow-out $G\cdot R\subseteq M$ and $G\times_{G_\nu}R$.  The
    contact form $\alpha$ on the flow-out can be reconstructed from
    the cross-section and the embedding $\iota:G_\nu\hookrightarrow
    G$.
  \end{enumerate}
\end{theorem}

\begin{remark}
  Note that the action of $G_\nu$ on the cross-section is in general
  not effective (even if the $G$-action on $M$ was).
\end{remark}

\begin{remark}
  The theorem uses the embedding $G_\nu\hookrightarrow G$.  If one
  considers a cross-section $R$ as an abstract $H$-manifold with
  $H\cong G_\nu$ and one embeds $H$ in two different ways into $G$
  ($\iota_1,\iota_2:\,H\hookrightarrow G$), then in general
  $G\times_{\iota_1 H}R\not\cong G\times_{\iota_2 H}R$. In the case of
  $\SO(3)$-manifolds however, the embedding of $\S^1$ into $\SO(3)$ is
  unique up to conjugation, and no problem will arise at this point.
\end{remark}

In the following corollary, the cross-section theorem will be applied
to $5$-dimensional contact $\SO(3)$-manifolds.

\begin{coro}
  Let $(M,\alpha)$ be a $5$-dimensional contact $\SO(3)$-manifold with
  moment map $\mu:\,M\to\so{3}^*$.  The cross-section $R$ is a
  $3$-dimensional contact $\S^1$-manifold without Legendrian orbits or
  fixed points.
  
  Conversely, let $(R,\alpha)$ be a $3$-dimensional contact
  $\S^1$-manifold without Legendrian orbits, and fixed points.  Then
  there is a $5$-dimensional contact $\SO(3)$-manifold $M$ that has
  $R$ as its cross-section.
\end{coro}
\begin{proof}
  The first part of the statement is a direct consequence of the
  cross-section theorem and Example~\ref{SO(3)-struktur}.  If $R$ had
  Legendrian orbits or fixed points, then $0$ would be contained in
  the image $\mu(R)$.
  
  For the second part, the manifold $M$ is given by
  $\SO(3)\times_{\S^1}R$, with the standard $\SO(3)$-action on the
  left factor.  The contact form on $M$ is constructed by taking
  $\alpha+\alpha(Z_R)\cdot Z^*$ on $\{e\}\times_{\S^1}R$, and moving
  it with the $\SO(3)$-action to the rest of $M$.  With $Z^*$, we mean
  the dual of $Z$ with respect to the standard basis $\{X,Y,Z\}$ of
  $\so{3}$.
\end{proof}

\begin{lemma}
  Let $(M,\alpha)$ and $(M^\prime,\alpha^\prime)$ be $5$-dimensional
  contact $\SO(3)$-manifolds.  An $\SO(3)$-contacto\-morphism
  $\Phi:\,M\to M^\prime$ induces an $\S^1$-contactomorphism between
  the cross-sections $R$ and $R^\prime$.
\end{lemma}
\begin{proof}
  The pull-back $\Phi^*\alpha^\prime$ is equal to $f\,\alpha$ with a
  positive function $f:\,M\to\R$.  For the moment maps, this gives
  $\mu^\prime\circ\Phi = f\cdot\mu$.  The restriction of $\Phi$ to $R$
  is then an $\S^1$-contactomorphism to $R^\prime$.
\end{proof}

\begin{lemma}\label{diffeo von R uebertraegt sich auf M}
  Let $(M,\alpha)$ and $(M^\prime,\alpha^\prime)$ be $5$-dimensional
  contact $\SO(3)$-manifolds, and let $R$ and $R^\prime$ be their
  respective cross-sections.  An $\S^1$-contactomorphism $\Phi:\,R\to
  R^\prime$ induces an $\SO(3)$-contactomorphism between the flow-outs
  $\SO(3)\cdot R\subset M$ and $\SO(3)\cdot R^\prime \subset
  M^\prime$.
\end{lemma}
\begin{proof}
  The map is given by
  $$
  \SO(3)\times_{\S^1}R\to \SO(3)\times_{\S^1}R^\prime,\quad [g,p]\mapsto
  [g,\Phi(p)]\;. $$
  One easily checks that these maps are well-defined, and
  respect the contact structures.
\end{proof}

\section{$5$-dimensional contact $\SO(3)$-manifolds}

The classification of closed symplectic $4$-manifolds with a
Hamiltonian $\SO(3)$- or $\SU(2)$-action was given in \cite{Iglesias}
and \cite{Audin}.  In the rest of the article, a proof to the theorem
below will given, which describes the classification of
$5$-dimensional contact $\SO(3)$-manifolds.

\begin{theorem}\label{mein hauptsatz}
  The following list gives a complete set of invariants for cooriented
  $5$-dimensional closed contact $\SO(3)$-manifold $M$, in the sense
  that there is an $\SO(3)$-contactomorphism between any two manifolds
  with equal invariants, and there exists a manifold for every choice
  of invariants from the list.
  \begin{itemize}
  \item The principal stabilizer is isomorphic to $\Z_k$ for some
    $k\in\N$ (including the trivial group, for $k=1$).
  \item The closure $\overline{R}$ of the cross-section is a compact
    $3$-dimensional contact $\S^1$-manifold without any fixed points
    or special exceptional orbits.  Each boundary component of
    $\overline{R}$ corresponds to a component of $\singorb{M}$. The
    orbits in the boundary are the only Legendrian orbits.
  \item If $M$ has singular orbits, then the principal stabilizer is
    either isomorphic to $\Z_2$ or trivial. In the first case, all
    components of $\singorb{M}$ are isomorphic to $\S^1\times\RP{2}$.
    If the principal stabilizer is trivial, one has two different
    types of components in $\singorb{M}$, which are either copies of
    $\S^1\times\S^2$ or $\S^1\widetilde\times\S^2 :=
    \R\times\S^2/\!\sim$ with the equivalence $(t,p)\sim (t+1,-p)$.  The
    Dehn-Euler number $n(R)$ is an integer, which describes how
    $\singorb{M}$ is glued onto $\regorb{M}$.  This Dehn-Euler number
    satisfies certain arithmetic conditions described in the
    Definition on page~\pageref{definition of dehn euler}.
  \end{itemize}
\end{theorem}

\begin{remark}\label{klassifikation von 3-mfkt}
  Contact $3$-dimensional $\S^1$-manifolds have been classified in
  \cite{Kamishima}.  The cross-section $R$ is thus determined by the
  following invariants:
  \begin{itemize}
  \item If $R$ is closed, it is determined solely by the genus of its
    orbit space $B:=R/\S^1$, the exceptional orbits, and the orbifold
    Euler number which cannot be zero.
  \item If $R$ is an open manifold, it is determined by the number of
    boundary components, the genus of its orbit space $B$, and its
    exceptional orbits.
  \end{itemize}
\end{remark}

Let $(M,\alpha)$ be a contact $5$-manifold and let $\SO(3)$ act by
contact transformations with moment map $\mu$.

\begin{lemma}\label{Hauptstabilisator}
  The principal stabilizer of a contact $\SO(3)$-manifold is
  isomorphic to $\Z_k$ for some $k\in\N$ (including the trivial group,
  for $k=1$).
\end{lemma}
\begin{proof}
  Since the moment map $\mu$ corresponding to the action is
  equivariant, $\stab(p)\le\mu(\stab(p))$.  The $\SO(3)$-structure of
  $\so{3}^*$ was given in Example~\ref{SO(3)-struktur}, and it follows
  that $\mu\equiv 0$ if the principal stabilizer is not one of $\Z_k$
  or $\S^1$.  But $\mu\equiv 0$ means that the action is trivial,
  which in particular contradicts effectiveness.
  
  In fact, the circle $\S^1$ can also be excluded: Assume $\exp(tX)$
  (for some $X\in\so{3}$, $X\ne 0$) leaves $p$ fixed, i.e.
  $\exp(tX)\cdot p=p$, then we have $\mu(p)=\mu(\exp(tX)\cdot
  p)=\Ad(\exp(-tX))^*\mu(p)$ and as a consequence $\ad(X)^*\mu(p) =
  0$.  Let now $X,Y,Z\in\so{3}$ be a standard basis of the Lie
  algebra.  Then, $\pairing{\mu(p)}{Z} = \pairing{\mu(p)}{[X,Y]} = 0$,
  $\pairing{\mu(p)}{Y} = -\pairing{\mu(p)}{[X,Z]} = 0$ and obviously
  $\pairing{\mu(p)}{X} = \alpha(X_M(p)) = 0$, i.e. $\mu(p) = 0$.
  
  Not only does this show that $\S^1$ cannot be a principal
  stabilizer, it also proves that all singular orbits lie in
  $\mu^{-1}(0)$, and the cross-section has no fixed points.
\end{proof}

The principal cross-section $R=\mu^{-1}(\R^+ Z^*)$ is a contact
$3$-manifold with a Hamiltonian $\S^1$-action.  The $\S^1$-orbits are
neither fixed points nor tangent to the contact structure. If
$0\notin\mu(M)$ the cross-section $R$ is a closed subset of $M$,
because $\R^+ Z^*\cap\mu(M)$ is compact, and hence $R$ is a closed
manifold and then $M$, as flow-out of $R$, is completely determined by
$R$.

\begin{lemma}\label{singorb ist nullmenge von mu}
  Let $(M,\alpha)$ be a $5$-dimensional contact $\SO(3)$-manifold.
  Then $\singorb{M}=\mu^{-1}(0)$.
\end{lemma}
\begin{proof}
  The preimage $\mu^{-1}(0)$ is the union of $\SO(3)$-orbits tangent
  to $\ker\alpha$, i.e.\ a collection of isotropic submanifolds. But
  isotropic submanifolds of a $5$-dimensional contact manifold have at
  most dimension $2$, and hence these orbits have to be singular.  On
  the other hand, the proof of Lemma~\ref{Hauptstabilisator} shows
  that all singular orbits lie in $\mu^{-1}(0)$.
\end{proof}

Furthermore a stabilizer of an exceptional orbit is isomorphic to some
$\Z_m$ and these orbits lie discrete surrounded by principal orbits.

\subsection{Examples}

In this section a few examples will be introduced that are continued
later in the article, while the theory is developed.

\begin{example}\label{sphaeren beispiel}
  The standard contact sructure on the $5$-sphere $\S^5\subset\C^3$
  is given at a point $(z_1,z_2,z_3)$ by
  $$
  \alpha_+ = \sum_{j=1}^3 \bigl(x_j\,dy_j-y_j\,dx_j\bigr)\;,
  $$
  with $z_j=x_j+iy_j$.  This contact form is invariant under the
  $\SO(3)$-action induced by the standard matrix representation.
  
  The stabilizer of a point $\mathbf{x}+i\mathbf{y}\in\S^5$ with
  $\mathbf{x}=(x_1,x_2,x_3)$ and $\mathbf{y}=(y_1,y_2,y_3)$ is the
  intersection of the stabilizer of $\mathbf{x}$ and that of
  $\mathbf{y}$.  If $\mathbf{x}$ and $\mathbf{y}$ are linearly
  independent, we have $\stab(\mathbf{x}+i\mathbf{y})=\{e\}$ and
  $\stab(\mathbf{x}+i\mathbf{y})\cong\S^1$ otherwise.
  
  For any matrix $A\in\so{3}$, the moment map is given by
  $\pairing{\mu(\mathbf{x}+i\mathbf{y})}{A} = 2\mathbf{x}^t
  A\mathbf{y}$.  The cross-section is then the set
  $$
  R=\{\mathbf{x}+i\mathbf{y}\in\S^5|x_1y_3-y_1x_3=x_2y_3-y_2x_3=0
  \text{ and } x_1y_2-y_1x_2>0\}\;.
  $$
  The condition $x_1y_2-y_1x_2>0$ implies that the other two
  equations, regarded as a linear system in $(x_3,y_3)$, have the
  unique solution $(x_3,y_3)=0$.  Hence the cross-section is given by
  $$
  R=\{(z_1,z_2,0)\in\S^5| x_1y_2-y_1x_2>0\}\;.
  $$
  
  The $\S^1$-action on $R$ is given by simultaneous rotations in the
  $(x_1,x_2)$- and $(y_1,y_2)$-plane.  Its orbit space $R/\S^1$ lies
  in a natural way in $\CP{1}$ with the projection $\pi:\, R\to
  R/\S^1$ given by $\pi(x_1+iy_1,x_2+iy_2,0)=[x_1+ix_2:y_1+iy_2]$.
  Note that the equation $x_1y_2-x_2y_1=0$ is well-defined in $\CP{1}$
  and its solutions are given by the standard embedding of $\RP{1}$.
  Hence $R/\S^1$ is diffeomorphic to an open disc and
  $R\cong\OpenDisc{2}\times\S^1$.

  Another $\SO(3)$-invariant contact form on $\S^5$ can be given by
  \begin{align*}
    \alpha_- &= i\,\sum_{j=1}^3 \bigl(z_j\,d\bar z_j - \bar
    z_j\,dz_j\bigr) \\
    & \qquad- i \bigl((z_1^2 + z_2^2 + z_3^2)\,d(\bar z_1^2 + \bar
    z_2^2 + \bar z_3^2) - (\bar z_1^2 + \bar z_2^2 + \bar
    z_3^2)\,d(z_1^2 + z_2^2 + z_3^2)\bigr)\;.
  \end{align*}
  Note that the first part of the form is identical to the standard
  form $\alpha_+$. It is easy to check that the second term does not
  give any contribution to the moment map, and hence $\mu_+ = \mu_-$.
  The cross-section for $\alpha_+$ and $\alpha_-$ are then of course
  also equal.

  The example will be continued at the end of the next section.
\end{example}

I would like to thank Otto van Koert for pointing out the following
examples to me.  As we will see later, these are all the simply
connected contact $\SO(3)$-manifolds with singular orbits of
dimension~$5$.  A good reference is \cite{Hirzebruch} and
\cite{Lutz_Meckert}.  The open book decomposition of these examples is
closely related to the $\SO(3)$-symmetry (\cite{KoertVan}).

\begin{example}\label{brieskorn beispiel}
  The Brieskorn manifolds $W^5_k\subset\C^4$ (with $k\in\N_0$) are
  defined as the intersection of the $7$-sphere with the zero set of
  the polynomial $f(z_0,z_1,z_2,z_3)= z_0^k+z_1^2+z_2^2+z_3^2$.  To
  make computations easier, assume the radius of the $7$-sphere to be
  $\sqrt{2}$.  It is well-known that $W^5_k$ is diffeomorphic to
  $\S^5$ for $k$ odd, and to $\S^2\times\S^3$ for $k$ even.
  
  Let $\SO(3)$ act linearly on $\C^4$ by leaving the first coordinate
  of $(z_0,z_1,z_2,z_3)$ fixed and multiplying the last three
  coordinates with $\SO(3)$ in its real standard representation, i.e.
  $A\cdot(z_0,z_1,z_2,z_3) := (z_0,A\cdot(z_1,z_2,z_3))$.  It is easy
  to check that this action restricts to $W^5_k$, because the
  polynomial $f$ can be written as $z_0^k + \norm{\mathbf{x}}^2 -
  \norm{\mathbf{y}}^2+2i\pairing{\mathbf{x}}{\mathbf{y}}$ with
  $\mathbf{x}=(x_1,x_2,x_3)$ and $\mathbf{y}=(y_1,y_2,y_3)$.  The only
  stabilizers that occur are $\{e\}$ and $\S^1$.  A point lies on a
  principal orbit, if and only if $\mathbf{x}$ and $\mathbf{y}$ are
  linearly independent.
  
  Finally the invariant $1$-forms
  \begin{align*}
    \alpha_k &= (k+1)\cdot(x_0\,dy_0-y_0\,dx_0) + 2\sum_{j=1}^3
    \left(x_j\,dy_j-y_j\,dx_j\right)
    \intertext{and}
    \alpha_{-k} &= -(k+1)\cdot(x_0\,dy_0-y_0\,dx_0) + 2\sum_{j=1}^3
    \left(x_j\,dy_j-y_j\,dx_j\right)
  \end{align*}
  are both of contact type on $W^5_k$.
  
  The infinitesimal generators of the $\SO(3)$-action do not have a
  $z_0$-component.  Hence the moment maps $\mu_k(z_0,z_1,z_2,z_3)$ for
  both $\alpha_k$ and $\alpha_{-k}$ are equal. They are given by
  $$
  \pairing{\mu_k}{X}=4(x_3y_2-x_2y_3),\quad\pairing{\mu_k}{Y}=4(x_1y_3-x_3y_1),
  \text{ and }\pairing{\mu_k}{Z}=4(x_2y_1-x_1y_2)\;.
  $$
  It can be seen with a similar computation as in
  Example~\ref{sphaeren beispiel} that the cross-section $R$ is given
  by the points $(z_0,z_1,z_2,0)\in W^5_k$ with $x_2y_1-x_1y_2>0$.
  
  The map $(z_0,z_1,z_2,0)\mapsto z_0$ from $R$ to the open unit disc
  is the projection of $R$ onto its quotient space (see
  \cite{Hirzebruch}).  The cross-section is $\S^1$-diffeomorphic to
  $\OpenDisc{2}\times\S^1$.

  The example will be continued at the end of the next section.
\end{example}

\subsection{Singular orbits}

In this section, we will show that each component of $\singorb{M}$
corresponds to one of three possible models.

\begin{lemma}\label{form fuer singulaere bahnen}
  Let $(M,\alpha)$ be a $5$-dimensional closed contact
  $\SO(3)$-manifold.  Recall from Lemma~\ref{Hauptstabilisator} that
  the principal stabilizer $H$ is either trivial or isomorphic to
  $\Z_k$.

  If $H\cong\Z_k$ with $k\ge 3$, then $M$ has no singular orbits.
  
  If $H\cong\Z_2$, then any component of $\singorb{M}$ has a
  neighborhood that is $\SO(3)$-diffeomorphic to a neighborhood of the
  zero-section in $\S^1\times T\RP{2}$, with trivial action on the
  first part and natural action on the second one.
  
  If $H$ is trivial, any component of $\singorb{M}$ has a neighborhood
  that is $\SO(3)$-diffeomorphic to a neighborhood of the zero-section
  in the vertical bundle $V\Etriv$ or $V\Etwist$, where $\Etriv$ is
  the trivial $\S^2$-bundle over $\S^1$ and $\Etwist$ is the twisted
  $\S^2$-bundle over $\S^1$.
  
  In all of these cases, there is up to $\SO(3)$-contactomorphisms a
  unique invariant contact form on sufficiently small neighborhoods of
  $\singorb{M}$.
\end{lemma}
In the rest of this section we will describe all possible cases, and
show the claims of the lemma.

One of the conclusion will be that the closure of the cross-section of
a $5$-dimensional contact $\SO(3)$-manifold $M$ is a compact
$3$-dimensional contact $\S^1$-manifold with boundary.  The interior
points of $R$ lie in regular $\SO(3)$-orbits, while $\partial R$ lies
in $\singorb{M}$.  The $\S^1$-orbits at the boundary are Legendrian.

\begin{lemma}[Equivariant Weinstein Theorem]\label{normalform fuer legendrebahn}
  Let $\orb(p)\hookrightarrow M$ be a Legendrian $\SO(3)$-orbit.  Then
  a neighborhood of $\orb(p)$ is $\SO(3)$-contactomorphic to a
  neighborhood of the zero-section in $(\R\oplus T^*\orb(p),
  dt+\lcan)$, where $\SO(3)$ acts by $g\cdot (t,v) = (t,g^{-1}_*v)$.
\end{lemma}
\begin{proof}
  There is an $\SO(3)$-invariant almost complex structure $J$ on the
  contact structure $\xi=\ker\alpha$ such that
  $$
  T_q\orb(p)\cap J\cdot (T_q\orb(p)) = \{0\} \text{ for all $q\in\orb(p)$.}
  $$
  The trivial line bundle $\epsilon^1$ spanned by the Reeb vector
  field of $\alpha$ is also $\SO(3)$-invariant.  This implies that the
  normal bundle of $T\orb(p)$ in $M$ can be equivariantly identified
  with $\epsilon^1\oplus T\orb(p)\cong\epsilon^1\oplus T^*\orb(p)$.
  The contact form restricts to $dt+c\,\lcan$ on the zero-section, and
  rescaling the fibre gives the desired form $dt+\lcan$.  This allows
  us to apply \cite[Theorem~5.2]{Lerman_Willett}, which states that
  there is a neighborhood of the orbit $\SO(3)$-contactomorphic to the
  normal bundle.
\end{proof}

By looking at the different stabilizers that can occur, it will be
seen that all singular orbits are either isomorphic to $\S^2$ with
stabilizer $\S^1$ or to $\RP{2}$ with stabilizer $\orthgroup(2)$.

\subsubsection{Fixed points}
The irreducible representations of $\SO(3)$ are all odd-dimensional.
This implies that $5$-dimensional contact $\SO(3)$-manifolds do not
have fixed points by the following argument. The vector space spanned
by the Reeb field is a trivial submodule of $T_pM$, and the contact
plane $(\xi_p,J_p)$ is a complex $2$-dimensional $\SO(3)$-module,
which also has to be trivial. That means the action on $T_pM$ is
trivial, which contradicts effectiveness.

\subsubsection{Stabilizer $\orthgroup(2)$}\label{singulaere bahnen mit O(2)}
The neighborhood of an orbit with stabilizer $\orthgroup(2)$ is
$\SO(3)$-equivariant to $\R\times T^*\orb(p)$ with $\orb(p)\cong
\RP{2}$.  The stabilizer of any non-zero element in $T^*\RP{2}$ is
isomorphic to $\Z_2$, which is then the principal stabilizer.

A connected component of $M_{(\orthgroup(2))}$ is an $\RP{2}$-bundle
over $\S^1$ (the closure $\overline{M_{(\orthgroup(2))}}$ is a closed
submanifold, possibly containing points with larger stabilizer than
$O(2)$, but we proved that $M$ has no fixed points, and hence
$\overline{M_{(\orthgroup(2))}} = M_{(\orthgroup(2))}$).  The
structure group of a $(G/H)$-bundle with the standard $G$-action on
the fibers are just the $G$-equivariant diffeomorphisms from $G/H$ to
itself.  It is not very difficult to see that these are given by
$N(H)/H$ (see~\cite{Bredon}).  In our case
$N(\orthgroup(2))/\orthgroup(2)=\orthgroup(2)/\orthgroup(2)=\{e\}$,
and hence every component of $M_{(\orthgroup(2))}$ is of the form
$\S^1\times\RP{2}$.  The neighborhood of such a component is
$\SO(3)$-diffeomorphic to $\S^1\times T^*\RP{2}$ with the standard
$\SO(3)$-action on the second part. A possible invariant contact form
is given by $dt+\lcan$, where $\lcan$ is the canonical $1$-form on
$T^*\RP{2}$.

In fact, the contact form above is the only one in a small
neighborhood of the singular orbit up to $\SO(3)$-contactomorphisms.
This can be proved in a similar way as Lemma~\ref{normalform fuer
  legendrebahn}: After pulling back the form to $\S^1\times
T^*\RP{2}$, one has $\alpha = f(t)\,dt + r(t)\,\lcan$ on the singular
orbits.  One can divide by $f(t)$ and then rescale the fibres to
obtain the standard form $dt + \lcan$, which allows us to use again
the Theorem from \cite{Lerman_Willett}.

In \ref{aequivalente mannigfaltigkeiten} and \ref{konstruktion von
  mannigfaltigkeiten}, it will be important to know how the
cross-section looks like in a neighborhood of the singular orbits.  We
compute the cross-section close to $M_{(\orthgroup(2))}$ in a
coordinate description.

A chart of $\RP{2}$ around $[1:0:0]$ is given by
$\R^2\to\RP{2},\,(q_1,q_2)\mapsto [1:q_1:q_2]$, and the
$\SO(3)$-action is induced by the standard matrix representation. Let
$X,Y,Z$ be the standard basis of $\so{3}$, where each element
generates the rotation around the corresponding axis of $\R^3$.  For
$Y$, for example the action looks like
\begin{align*}
  \exp(t Y)\cdot[1:q_1:q_2] &= [\cos t + q_2\sin t : q_1:q_2\cos t - \sin t]\\
  &= \left[1:\frac{q_1}{\cos t + q_2\sin t}:\frac{q_2\cos t - \sin t}{\cos t + q_2\sin t}\right]
\end{align*}
The infinitesimal generators of the action are given in this chart by
\begin{align*}
  X_{\RP{2}}([1:q_1:q_2]) &= q_2\,\partial_{q_1}-q_1\,\partial_{q_2}\;, \\
  Y_{\RP{2}}([1:q_1:q_2]) &= -q_1q_2\,\partial_{q_1}-(1+q_2^2)\,\partial_{q_2}\;, \\
  Z_{\RP{2}}([1:q_1:q_2]) &=
  -(1+q_1^2)\,\partial_{q_1}-q_1q_2\,\partial_{q_2}\;.
\end{align*}
and the moment map is
\begin{align*}
  \pairing{\mu(t,q_1,q_2,p_1,p_2)}{X} &= q_2p_1 - q_1p_2 \;, \\
  \pairing{\mu(t,q_1,q_2,p_1,p_2)}{Y} &= -q_1q_2p_1 - (1+q_2^2)\,p_2 \;, \\
  \pairing{\mu(t,q_1,q_2,p_1,p_2)}{Z} &= -(1+q_1^2)\,p_1 - q_1q_2p_2 \;.
\end{align*}
Elements of $\mu^{-1}(\R^+Z^*)$ have $p_1\ne 0$ or $p_2\ne 0$, and for
such elements $q_2p_1 - q_1p_2 = 0$ and $-q_1q_2p_1 - (1+q_2^2)\,p_2 =
0$ hold.  These two equations can be read as a linear system in $p_1$
and $p_2$, and there are only non-trivial solutions if the
corresponding determinant vanishes, that is, if $-q_2\,(1+q_2^2) -
q_1^2q_2 = -q_2\,(1+q_1^2+q_2^2) = 0$.  If this is the case, then
$q_2=0$, and from this it follows that $p_2=0$.  The cross-section $R$
consists of vectors in $T\RP{2}$ tangent to $\RP{1}$, but pointing
only in positive direction (with the embedding of $\RP{1}$ in $\RP{2}$
given by $[a:b]\mapsto[a:b:0]$).

The restriction of the contact form on $R$ is given in the chart above
by $dt+p_1\,dq_1 $.  Hence $\alpha$ is of contact type even on the
boundary of $\overline{R}$, and the orbits of the $\S^1$-action are
Legendrian on $\partial\overline{R}\cong\S^1\times\S^1$.

A collar neighborhood of $\partial\overline{R}$ is of the form
$\S^1\times[0,\epsilon)\times\S^1$ with contact form $dt+r\,d\phi$ and
action $e^{i\theta}\cdot(t,r,\phi)=(t,r,\phi+2\theta)$. The embedding
of this neighborhood into $M$ is given by
$$
(t,r,\phi)\mapsto \big(t,[\cos(\phi/2):\sin(\phi/2):0],
-r\sin(\phi/2)\,\partial_1 + r\cos(\phi/2)\,\partial_2\big)\;,
$$
and the points $(t,0,0)\in\partial\overline{R}$ all have equal
stabilizer in $\SO(3)$.

\subsubsection{Stabilizer $\S^1$}\label{singulaere bahnen mit SO(2)}
The neighborhood of such an orbit is $\SO(3)$-diffeomorphic to
$\R\times T\S^2$ with trivial action on the first and standard action
on the second component.  The principal stabilizer is trivial.  A
connected component of $M_{(\SO(2))}$ is a closed manifold, because no
fixed points or points with stabilizer $\orthgroup(2)$ do exist, and
hence $M_{(\SO(2))}$ is diffeomorphic to an $\S^2$-bundle over $\S^1$.
The structure group of such a bundle is $N(\SO(2))/\SO(2)\cong\Z_2$,
hence the only two $\S^2$-bundles over $\S^1$ are the trivial one
$\Etriv$ and the twisted one $\Etwist$.  They can be described by the
equivalence relations $(t,p)\sim (t+1,p)$ and $(t,p)\sim (t+1,-p)$
(with $t\in\R$ and $p\in\S^2$) respectively.  A neighborhood of a
component of $\singorb{M}$ is diffeomorphic to the corresponding
vertical bundle.  The $\SO(3)$-action on the second component of
$\R\times\S^2$ is compatible with these identifications, and one
obtains an action on either vertical bundle $V\Etriv$ and $V\Etwist$.

A possible invariant contact form is given by $dt+\lcan$ on $\R\times
T^*\S^2$, where $T^*\S^2$ is identified with $T\S^2$ via an invariant
metric.  This form descends to $V\Etriv$ and also to $V\Etwist$,
because the reflection in the construction of $\Etwist$ is induced by
a diffeomorphism of $\S^2$, and $\lcan$ on $T^*N$ remains invariant
under maps induced by diffeomorphisms of the base space $N$.

In a small neighborhood of $M_{(\SO(2))}$, every invariant contact
form is $\SO(3)$-contactomorphic to $dt+\lcan$.  The proof of this
fact is completely analogous to the one for orbits with stabilizer
$\orthgroup(2)$ above, and will be omitted.

Now we will describe how the cross-section looks like in a
neighborhood of the singular orbits. The moment map $\mu$ is given in
the neighborhood of a singular orbit by
$$
\pairing{\mu(t,q,p)}{X} = p^t Xq
$$
with $(t,q,p)\in\R^1\times T^*\S^2 \subseteq
\R^1\times\R^3\times\R^3$ and $X\in\so{3}$ in its standard matrix
representation.  One easily checks that the cross-section is the set
of points $(t,q,p)$ where $q$ lies in the equator of the sphere and
$p$ is a vector tangent to the equator at $q$, with all these vectors
oriented the same way.  The $\S^1$-action on the cross-section is
induced by rotations around the $z$-axis of the sphere.

\begin{figure}[htbp]
  \begin{picture}(0,0)%
    \includegraphics{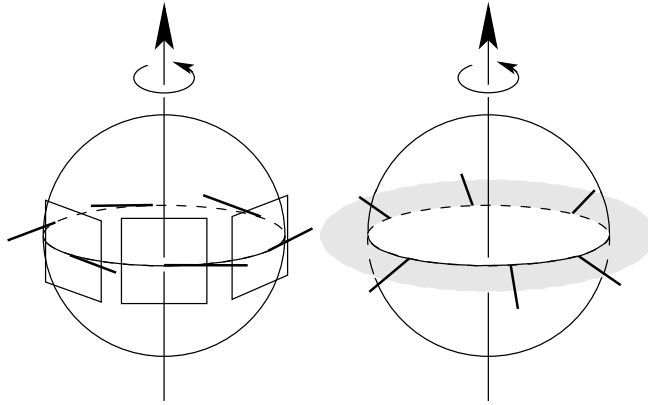}%
  \end{picture}%
  \setlength{\unitlength}{4144sp}%
  \begingroup\makeatletter\ifx\SetFigFont\undefined%
  \gdef\SetFigFont#1#2#3#4#5{%
    \reset@font\fontsize{#1}{#2pt}%
    \fontfamily{#3}\fontseries{#4}\fontshape{#5}%
    \selectfont}%
  \fi\endgroup%
  \begin{picture}(3914,2418)(-1521,-1098)
  \end{picture}

  \caption{On the left the cross-section around an exceptional orbit is
    displayed: It consists of vectors at the equator pointing into
    positive direction. The picture on the right displays a model more
    accessible to the imagination: The cross-section sits as a ring
    around the equator of the sphere. Vectors pointing into the
    cross-section are normal to the sphere.}
\end{figure}

For $\Etriv$, a collar neighborhood of the boundary
$\partial\overline{R}$ can be given by
$\S^1\times[0,\epsilon)\times\S^1$, while for components of type
$\Etwist$, the form $\R\times[0,\epsilon)\times\S^1/\!\sim$ with the
equivalence relation $(t,r,\phi)\sim(t+1,r,\phi+\pi)$ will be used.
The contact form is $dt+r\,d\phi$ in both cases, and the $\S^1$-action
is $e^{i\theta}\cdot(t,r,\phi)=(t,r,\phi+\theta)$.  The embedding of
$\overline{R}$ into the neighborhood of $\singorb{M}$ is given by
$$
(t,r,\phi)\mapsto \big(t; (\cos\phi,\sin\phi,0);r\cdot(-\sin\phi,\cos\phi\,0)\big)\;.
$$
With this embedding, the points $(t,0,0)$ and $(t,0,\pi)$ in
$\partial\overline{R}$ all have equal stabilizer.

This concludes the description of all singular orbits, and the proof
of Lemma~\ref{form fuer singulaere bahnen}.

\begin{sphaerenbeispiel}[cont.]
  As described above, the singular orbits of $\S^5$ are composed of
  all points $\mathbf{x}+i\mathbf{y}$ where $\mathbf{x}=(x_1,x_2,x_3)$
  and $\mathbf{y}=(y_1,y_2,y_3)$ are linearly dependent. The singular
  orbits are $2$-spheres, and we have to decide whether the component
  of $\singorb{\S^5}$ is equal to $\Etriv$ or to $\Etwist$. This of
  course is independent of the contact structure. The only points
  invariant under rotations around the $z_3$-axis are
  $(0,0,e^{i\phi})$ with $0\le\phi<2\pi$. But since $(0,0,1)$ and
  $(0,0,-1)$ both lie in $\orb(0,0,1)$, we have
  $\singorb{\S^5}\cong\Etwist$.
\end{sphaerenbeispiel}

\begin{brieskornbeispiel}[cont.]
  Now we will determine the type of the singular orbits of $W_k^5$.
  This of course does not depend on the contact structure.  As we said
  above, a point $(z_0,z_1,z_2,z_3)\in W_k^5$ lies on a singular orbit
  if and only if $\mathbf{x}$ is parallel to $\mathbf{y}$, where
  $\mathbf{x}=(x_1,x_2,x_3)$ and $\mathbf{y}=(y_1,y_2,y_3)$.  In
  particular, consider the points that are invariant under rotations
  around the $z_1$-axis.  They are given by $\bigl\{\bigl(e^{i\phi},
  \pm ie^{\frac{ki}{2}\phi},0,0\bigr) |\,0\le\phi<2\pi\bigr\}$. For
  $k$ odd, all points lie on a single path, but for $k$ even, there
  are two connected components.  Hence, one obtains
  $\singorb{\left(W_k^5\right)}\cong \Etwist$ for $k$ odd, and
  $\singorb{\left(W_k^5\right)}\cong \Etriv$ for $k$ even.
  
  So far all invariants found for $(W_k^5,\alpha_{\pm k})$, and
  $(W_{k^\prime}^5,\alpha_{\pm k^\prime})$ are equal if $k\equiv
  k^\prime\mod 2$.  But at the end of the next section, a last
  invariant will be computed that allows us to distinguish all of the
  $(W_k^5,\alpha_{\pm k})$.
\end{brieskornbeispiel}

\subsection{Equivalence between contact $\SO(3)$-manifolds}\label{aequivalente mannigfaltigkeiten}

In this section, the necessary and sufficient conditions for the
existence of an $\SO(3)$-equivariant contactomorphism $\Phi:\,M\to
M^\prime$ between two $5$-dimensional contact $\SO(3)$-manifolds
$(M,\alpha)$ and $(M^\prime,\alpha^\prime)$ will be given.

If there are no singular orbits on $M$, then $0\notin\mu(M)$ and the
whole manifold is determined according to Theorem~\ref{cross-section}
by its cross-section.  Two contact $5$-manifolds with an
$\SO(3)$-action without singular orbits are thus equivalent if and
only if their cross-sections are.  The possible cross-sections, being
closed contact $3$-manifold with $\S^1$-actions, have been classified
in \cite{Kamishima}.

On the other hand, if $0\in\mu(M)$, then
$M=\regorb{M}\cup\singorb{M}$, but there are several ways to glue both
parts.  The flow-out $\SO(3)\cdot R\cong \SO(3)\times_{\S^1}R$ is
determined by $R$, but for the whole of $M$ the problem is that
$p\in\partial\overline{R}$ does not ``remember'' as point in the
$\S^1$-manifold $\overline{R}$, which stabilizer $\stab(p)\le\SO(3)$
it had in $M$.

The solution lies in choosing an arbitrary point
$p_0\in\partial\overline{R}$ and marking all other points $p$ in the
boundary with $\stab(p)=\stab(p_0)\le\SO(3)$.  The marked points form
curves in $\partial\overline{R}$.  If the boundary component
corresponds to $\Etriv$, these curves are given by two sections to the
$\S^1$-action that are related to each other by a
$180^\circ$-rotation.  If the component corresponds to $\Etwist$, the
marked points lie on a single curve, which intersects each
$\S^1$-orbit twice.  If the singular orbits have stabilizer isomorphic
to $\orthgroup(2)$, then the marked points form a single section.

Another way to describe the situation is the following: Gluing
$\singorb{M}$ onto $\regorb{M}$ can be achieved by gluing $R$ onto the
cross-section in the neighborhood of $\singorb{M}$.  This means that
one has to identify two tori.  The generators of the homology in
$\partial R$ are given by an $\S^1$-orbit and a section $\sigma$ to
the $\S^1$-action in $\overline{R}$.  The generators of the homology
of $\overline{R}\cap\singorb{M}$ can be described by an $\S^1$-orbit,
and by a curve of marked points as fixed above.  The $\S^1$-obits have
to coincide in both parts, and the only freedom when gluing consists
in choosing the relative position of the other two homology classes.

\begin{lemma}\label{kontaktomorphismus auf R gibt kontaktomorphismus auf M}
  Let $(M,\alpha)$ and $(M^\prime,\alpha^\prime)$ be two
  $5$-dimensional contact $\SO(3)$-manifolds with principal
  cross-sections $(R,\alpha)$ and $(R^\prime,\alpha^\prime)$. Assume
  there is an $\S^1$-contactomorphism $\psi$ between $\overline{R}$
  and $\overline{R^\prime}$ that maps the marked curves
  $\gamma_1,\ldots,\gamma_n$ in $\partial\overline{R}$ onto the marked
  curves in $\partial\overline{R^\prime}$, i.e.\ 
  $\psi\circ\gamma_i=\gamma_i^\prime$. Then there is an
  $\SO(3)$-equivariant contactomorphism $\Psi:\,M\to M^\prime$.
\end{lemma}
\begin{proof}
  Over the flow-out $\SO(3)\cdot R$ and $\SO(3)\cdot R^\prime$ the
  claim holds. Hence if $\singorb{M}=\emptyset$, then the statement is
  true. The problem for $\partial\overline{R}\ne\emptyset$ is that
  $\psi$ extends to an $\SO(3)$-homeomorphism on $M$, but this map is
  in general not smooth at the singular orbits. Hence we will need to
  deform $\psi$ in a neighborhood of~$\partial\overline{R}$.
  
  Choose a component $K$ of $\singorb{M}$.  The image $\psi(K)$ in
  $\singorb{M^\prime}$ is of the same type: If the principal
  stabilizer of $R$ is isomorphic to $\Z_2$, then every component in
  $\singorb{M}$ and $\singorb{M^\prime}$ is diffeomorphic to
  $\S^1\times\RP{2}$, and if the principal stabilizer of $R$ is
  trivial, then the two types of component in $\singorb{M}$ and
  $\singorb{M^\prime}$ can be distinguished by the curves of marked
  points.

  Now one can represent the neighborhood of $K$ and $\psi(K)$ by the
  standard models described at the end of Section~\ref{singulaere
    bahnen mit O(2)} and \ref{singulaere bahnen mit SO(2)}.  The
  cross-section is either given by
  $(\R\times[0,c)\times\S^1/\!\sim,dt+r\,d\phi)$ for $\Etwist$ or by
  $(\S^1\times[0,c)\times\S^1,dt+r\,d\phi)$ for the other two types of
  singular orbits.

  The map $\psi$ is $\S^1$-equivariant, thus
  $$
  \psi(t,r,\phi) = \big(T(t,r), R(t,r), \phi+\Phi(t,r)\big)\;.
  $$
  Furthermore it rescales the form $\alpha=dt+r\,d\phi$ by a
  function $f(t,r) > 0$, i.e.\ 
  $$
  f(t,r)\,dt+r f(t,r)\,d\phi = f\alpha=\psi^*\alpha =
  \left(\frac{\partial T}{\partial t} +
    R\cdot\frac{\partial\Phi}{\partial t}\right)\,dt + R\,d\phi +
  \left(\frac{\partial T}{\partial
      r}+R\cdot\frac{\partial\Phi}{\partial r}\right)\, dr\;.
  $$
  The consequences are $R(t,r) = r f(t,r)$, $\partial_tT(t,r) + r
  f(t,r)\cdot\partial_t\Phi(t,r) = f(t,r)$, and $\partial_r T(t,r) + r
  f(t,r)\cdot\partial_r\Phi(t,r) = 0$.  The boundary is mapped onto
  the boundary, i.e.\ $R(t,0)=0$. We can assume $T(0,0) = 0$ and
  $\Phi(0,0) = 0$. Also, all of the three cases $\Etriv$, $\Etwist$,
  and $\S^1\times\RP{2}$ lead to $\Phi(t,0)=0$, because the $\gamma_i$
  are mapped onto the $\gamma_i^\prime$.

  Let $\rho_\epsilon:\R^+\to [0,1]$ be the smooth map
  \begin{equation*}
    \rho_\epsilon(r) =
    \begin{cases}
      0 & \text{for $r\le \epsilon/2$} \\
      N(\epsilon)\cdot\int_{\epsilon/2}^r \exp\frac{\epsilon^2}{4(x-\epsilon/2)(x-\epsilon)}\,dx & \text{for $\epsilon/2 < r < \epsilon$} \\
      1 & \text{for $r\ge\epsilon$}
    \end{cases}
  \end{equation*}
  with $N(\epsilon)$ the reciprocal value of
  $\int_{\epsilon/2}^\epsilon
  \exp\frac{\epsilon^2}{4(x-\epsilon/2)(x-\epsilon)}\,dx$.  The
  maximum of the derivative of this function is
  $N(\epsilon)\cdot\exp(-4) = N(1)e^{-4}/\epsilon$.  One can now
  replace the original map $\psi$ by
  $$
  \widehat\psi(t, r, \phi) := \bigl(T(t,r), R(t,r),
  \phi+\rho_\epsilon(r)\cdot\Phi(t,r)\bigr)\;.
  $$
  It is easy to check that $\widehat\psi$ is well-defined on the
  cross-section $R$: The relations $\psi(t+2\pi a,r,\phi + 2\pi b) =
  \psi(t,r,\phi)+(2\pi a, 0, 2\pi b)$ carry over to $\widehat\psi$.
  
  The map $\widehat\psi$ is equal to $\big(T(t,r), rf(t,r), \phi\big)$
  for points with $r\le\epsilon/2$ and equal to $\psi$ for points with
  $r\ge\epsilon$. It is also an $\S^1$-diffeomorphism. The determinant
  of the differential $d\widehat\psi$ is equal to the one of $d\psi$.
  The injectivity and surjectivity follow easily from the same
  properties of $\psi$.  For example to show that
  $(t^\prime,r^\prime,\phi^\prime)$ lies in the image of
  $\widehat\psi$, use that there is a $(t,r,\phi)$ with
  $\psi(t,r,\phi)=(t^\prime,r^\prime,\phi^\prime)$. Then
  $\widehat\psi\bigl(t,r,\phi +
  (1-\rho_\epsilon(r))\cdot\Phi(t,r)\bigr) =
  (t^\prime,r^\prime,\phi^\prime)$.
  
  There is now an $\SO(3)$-diffeomorphism $\widehat\Psi$ on $M$
  extending $\widehat\psi$. Away from the singular orbits, the map
  $\widehat\Psi$ is given as in the proof of Lemma~\ref{diffeo von R
    uebertraegt sich auf M}.  In the neighborhood of $\singorb{M}$ one
  can use the standard model for $\Etriv$ and $\Etwist$, where the map
  $\widehat\Psi$ is given by
  $$
  \widehat\Psi:\,(t;p,v)\mapsto (T(t,\norm{v});\,p,\, f(t,\norm{v})\,v)\;,
  $$
  for $p\in\S^2$ and for $v\in T^*_p\S^2$ with
  $\norm{v}<\epsilon/2$.  If the component of $\singorb{M}$ was
  diffeomorphic to $\S^1\times\RP{2}$ the map is given by the
  projectivization of $\widehat\Psi$ defined above.  These maps
  clearly define $\SO(3)$-equivariant diffeomorphisms in the
  neighborhood of a singular orbit, but one still needs to check that
  this definition is compatible with the map given in the proof of
  Lemma~\ref{diffeo von R uebertraegt sich auf M}.  Because both maps
  are $\SO(3)$-equivariant, it is enough to check that these maps
  agree on the cross-section $R$.  But $\widehat\Psi$ restricted to
  $R$ gives back the map $\widehat\psi$.  This shows that
  $\widehat\Psi$ is a globally defined map.
  
  The map $\widehat\Psi$ is an $\SO(3)$-diffeomorphism, but it is only
  a contactomorphism far away from the singular orbits. All of the
  $\SO(3)$-invariant $1$-forms in the family $\alpha_s :=
  (1-s)\,\alpha + s\,\widehat\Psi^*\alpha$ on $M$ satisfy the contact
  condition.  This can easily be checked in a small neighborhood of
  the singular orbits by using the local form given above.  On
  $\princorb{M}$, one checks the contact condition along $R$ (by
  choosing $\epsilon$ small enough) and then uses $\SO(3)$-invariance.
  The equivariant Gray stability shows that $\widehat\Psi$ deforms to
  an $\SO(3)$-contactomorphism~$\Psi$.
\end{proof}

Of course, the next question is how to find maps with the properties
required in Lemma~\ref{kontaktomorphismus auf R gibt
  kontaktomorphismus auf M}.  For this, we need to define a last
invariant for the cross-section.

Let $\overline{R}$ be a compact oriented $3$-dimensional
$\S^1$-manifold with non-empty boundary.  Denote the components of
$\partial\overline{R}$ by $\partial\overline{R}_j$ ($j=1,\ldots,N$)
and assume that on each of the boundary components a smooth closed
curve $\gamma_j$ is given that intersects the $\S^1$-orbits
transversely.  Orient the curves in such a way that $\dot\gamma_j$
followed by the inifinitesimal generator $Z_R$ of the $\S^1$-action
gives the orientation of $\partial\overline{R}_j$.

The $\gamma_j$ should be of the same form as the marked points
described above, i.e.\ if the principal stabilizer is isomorphic to
$\Z_2$, assume $\gamma_j$ intersects each $\S^1$-orbit in $\partial
\overline{R}_j$ exactly once.  If the principal stabilizer of $R$ is
trivial, the curves are either sections or intersect each orbit twice.

On the boundary of a small tubular neighborhood of the exceptional
orbits one can define standard sections (\cite{Orlik}), which can be
extended to a global section $\sigma$ of $\overline{R}\to
\overline{R}/\S^1$.  Let $\sigma$ be oriented in such a way that the
tangent space to the image of $\sigma$ followed the generator of the
$\S^1$-action gives the positive orientation of~$\overline{R}$.

\begin{defi}\label{definition of dehn euler}
  Denote the intersection number of two oriented loops $\alpha$ and
  $\beta$ in an oriented torus by $\intersection{\alpha}{\beta}$.  If
  the principal stabilizer in $R$ is trivial define the
  \textbf{Dehn-Euler-number}\index{Dehn-Euler number}
  $n(R,\gamma_1,\ldots,\gamma_N)\in\Z$ by
  \begin{align*}
    n(R,\gamma_1,\ldots,\gamma_N) &:= 2\sum_{j=0}^m
    \intersection{\gamma_j}{\partial\sigma} + \sum_{j=m+1}^N
    \intersection{\gamma_j}{\partial\sigma}\;,
  \end{align*}
  where we assume the first $m$ curves to be sections to the
  $\S^1$-action, and the other curves to intersect each orbit twice.
  Note that the first term is a sum over even numbers and the second
  term is a sum over odd numbers.
  
  If the principal stabilizer is isomorphic to $\Z_2$ define the
  Dehn-Euler number by
  \begin{align*}
    n(R,\gamma_1,\ldots,\gamma_N) &:= \sum_{j=1}^N
    \intersection{\gamma_j}{\partial\sigma}\;.
  \end{align*}
  In this case $n(R,\gamma_1,\ldots,\gamma_N)$ can be any integer.
\end{defi}

The Dehn-Euler number is very similar to the Euler invariant for an
$\S^1$-manifold.  To see that $n(R,\gamma_1,\ldots,\gamma_N)$ is
independent of the section chosen, assume two different sections
$\sigma_1$ and $\sigma_2$ (that are homotopic to the standard sections
around the exceptional orbits) are given.

There is a function $f:\,R/\S^1\to\S^1$, such that
$\sigma_2(p)=\sigma_1(p)\cdot f(p)$.  The \textbf{rotation number}
$\rot(\left.f\right|_{\partial R_j})$ is defined as the degree of the
map $\left.f\right|_{\partial R_j}:\,\partial R_j/\S^1\cong
\S^1\to\S^1$.  The sum $\sum \rot(\left.f\right|_{\partial R_j})$ over
all boundary components of $R$ vanishes, because the degree of a map
$\Disc{2}\to\S^1$ vanishes on $\partial\Disc{2}$. We can cut $R/\S^1$
open to obtain a disc, and the extra contributions from the cuts
cancel out.  With the equations
$\intersection{\gamma_j}{\partial\sigma_2}-\intersection{\gamma_j}{\partial\sigma_1}=\rot(\left.f\right|_{\partial
  R_j})$ for $j\le k$, and
$\intersection{\gamma_j}{\partial\sigma_2}-\intersection{\gamma_j}{\partial\sigma_1}=2\rot(\left.f\right|_{\partial
  R_j})$ for $k<j\le N$, it follows that
\begin{multline*}
  2\sum_{j=0}^k (\intersection{\gamma_j}{\partial\sigma_1}-\intersection{\gamma_j}{\partial\sigma_2}) +  \sum_{j=k+1}^N (\intersection{\gamma_j}{\partial\sigma_1} - \intersection{\gamma_j}{\partial\sigma_2}) \\
  = 2\sum_{i=0}^k \rot(f\left|_{\partial R_i}\right.) + \sum_{j=k+1}^N
  2\rot(f\left|_{\partial R_j}\right.) = 0
\end{multline*} 

Note also that the orientation of the $\S^1$-action has no effect on
$n(R,\gamma_1,\ldots,\gamma_N)$.  To compute
$n(R,\gamma_1,\ldots,\gamma_N)$ we can use again the section $\sigma$,
because the standard sections around the exceptional orbits do not
change with the orientation of the $\S^1$-action.  The direction of
the boundary curves $\gamma_j$ and the orientation of $\sigma$ are
inverted.  But then the intersection number remains unchanged.

\begin{remark}
  In Lemma~\ref{diffeo von R uebertraegt sich auf M}, it was shown
  that the cross-section $R$ (as contact $\S^1$-manifold) is an
  invariant of a $5$-dimensional contact manifold $M$.  It has just
  been proved that the number $n(R,\gamma_1,\ldots,\gamma_m)$ is also
  an invariant of $M$, because under an $\SO(3)$-contactomorphism the
  marked curves are mapped onto each other.  Below we will now finish
  the proof that a manifold $M$ is completely determined by the
  invariants mentioned in Theorem~\ref{mein hauptsatz} (i.e.\ 
  cross-section, singular orbits and $n(R)$).
\end{remark}

The $3$-manifolds in the following lemma are cross-sections of
$5$-manifolds.

\begin{lemma}
  Let $(R,\alpha)$ and $(R^\prime,\alpha^\prime)$ be two
  $\S^1$-diffeomorphic $3$-dimensional contact $\S^1$-manifolds
  without fixed points, but both with $N$ boundary components.  Let
  the orbits in the boundary be the only ones that are Legendrian.
  Assume further that on each of the boundary components $\partial
  R_j$ and $\partial R_i^\prime$ curves $\gamma_j$ and
  $\gamma_i^\prime$ are specified such that for both manifolds the
  first $k$ curves ($k\le N$) are sections to the $\S^1$-action and
  the other curves intersect each orbit exactly twice.  Then there is
  an $\S^1$-contactomorphism $\Phi:\,R\to R^\prime$ such that
  $\Phi\circ\gamma_j = \gamma_j^\prime$, if and only if
  $n(R,\gamma_1,\ldots,\gamma_N) =
  n(R^\prime,\gamma_1^\prime,\ldots,\gamma_N^\prime)$.
\end{lemma}
\begin{proof}
  The basic strategy is to find diffeomorphic sections with certain
  properties in $R$ and $R^\prime$. With these sections one can
  construct an $\S^1$-diffeomorphism between the $3$-manifolds that
  maps the boundary curves in $R$ onto the ones in $R^\prime$.
  Afterwards this map is deformed to obtain a contactomorphism.
  
  By \cite{Kamishima}, the contact form around an exceptional orbits
  is locally unique up to $\S^1$-contacto\-mor\-phisms.  Thus one can
  start the construction of $\Phi$ by taking an
  $\S^1$-contactomorphism from a small neighborhood of the exceptional
  orbits in $R$ to a neighborhood of the orbits of the same type in
  $R^\prime$.  Choose also, for each $j\in\{1,\ldots,N-1\}$, an
  $\S^1$-diffeomorphism from a neighborhood of $\partial R_j$ to a
  neighborhood of $\partial R_j^\prime$ that maps $\gamma_j$ onto
  $\gamma_j^\prime$ .
  
  The standard sections to the $\S^1$-action around the exceptional
  orbits extend to a global section $\sigma$ on $\princorb{R}$.  In
  $R^\prime$, construct a section in the following way: Take $\sigma$
  in the neighborhood of the exceptional orbits and in the
  neighborhood of $\partial R_j$ for $1\le j\le N-1$ and map it with
  $\Phi$ to $R^\prime$.  Now extend the image of $\sigma$ to a global
  section $\sigma^\prime$ on $\princorb{R^\prime}$.
  
  By the assumptions of the lemma, we know that
  $n(R,\gamma_1,\ldots,\gamma_N) =
  n(R^\prime,\gamma_1^\prime,\ldots,\gamma_N^\prime)$, and by our
  construction $\intersection{\sigma}{\gamma_j} =
  \intersection{\sigma^\prime}{\gamma_j^\prime}$ for all $1\le j \le
  N-1$.  It follows that the intersection numbers
  $\intersection{\sigma}{\gamma_N}$ and
  $\intersection{\sigma^\prime}{\gamma_N^\prime}$ are also equal.
  Hence one can homotope $\sigma^\prime$ in such a way that its
  position with respect to $\gamma_N^\prime$ is the same as the one of
  $\sigma$ with respect to $\gamma_N$.
  
  One can map $\sigma$ onto $\sigma^\prime$ and by using the
  $\S^1$-action, we obtain an $\S^1$-diffeomorphism $\Phi:\, R\to
  R^\prime$, such that $\Phi\circ\gamma_j=\gamma_j^\prime$ for all
  $j=1,\ldots,N$.
  
  To transform the map above into a contactomorphism we need to
  sharpen an argument given in \cite{Lutz} and \cite{Kamishima} to
  avoid moving the curves on the boundaries.  The neighborhoods of the
  boundaries are of the form $\S^1\times[0,\delta)\times\S^1$ with
  coordinates $(t,r,\phi)$, and the circle action on the last
  coordinate.  Assume one contact form to be $\alpha = dt+r\,d\phi$
  and the other one $\alpha^\prime = g(t,r)\,dt + h(t,r)\,dr +
  f(t,r)\,d\phi$.  The orbits in the boundary are Legendrian, hence
  $f(t,0)=0$ and $\partial_t f(t,0) = 0$.  Thus the contact condition
  along such an orbit becomes $g(t,0)\ne0$, and we can divide the
  whole form by the function $g$ to obtain the equivalent form $dt +
  h(t,r)\,dr + f(t,r)\,d\phi$ (with new functions $f$ and $h$).
  
  Define now a map $\Psi:\,R\to R$ by
  $$
  (t,r,\phi)\mapsto (t - (1-\rho_\epsilon(r))r h(t,0),r,\phi)
  $$
  for points with $r<\epsilon$ and the identity otherwise.  Here
  $\rho_\epsilon$ is the map defined in the proof of
  Lemma~\ref{kontaktomorphismus auf R gibt kontaktomorphismus auf M}.
  
  The map $\Psi$ is an $\S^1$-diffeomorphism.  It is surjective,
  because it is the identity on the two tori
  $\S^1\times\{0\}\times\S^1$ and $\S^1\times\{\epsilon\}\times\S^1$.
  The map is a local diffeomorphism because $\det(d\Psi) = 1 -
  r(1-\rho_\epsilon(r))\,\partial_t h(t,0)$ does not vanish, if we
  choose $\epsilon$ small enough.  Injectivity relies on a similar
  argument: If $\Psi(t,r,\phi) = \Psi(t^\prime,r^\prime,\phi^\prime)$,
  then clearly $\phi = \phi^\prime$ and $r = r^\prime$.  Finally
  $t-t^\prime = r (1-\rho_\epsilon(r))\, (h(t,0)-h(t^\prime,0))$.
  With the mean value theorem one sees that if $t \ne t^\prime$, one
  has $1 = r (1-\rho_\epsilon(r))\,\partial_t h(\hat t,0)$ with $\hat
  t\in(t,t^\prime)$, which is not possible if $\epsilon$ is chosen
  small enough.
  
  For $r=0$ the forms $\alpha$ and $\Psi^*\alpha^\prime$ are equal,
  hence the linear interpolation $\alpha_s = (1-s)\,\alpha +
  s\,\Psi^*\alpha^\prime$ consists of $\S^1$-invariant contact forms.
  To apply the Moser trick one considers the vector field $X_s$ that
  is the solution to the equations
  $$
  \iota_{X_s}\alpha_s \quad\text{ and }\quad \iota_{X_s}d\alpha_s =
  \lambda_s\,\alpha_s -\dot\alpha_s, $$
  with the function $\lambda_s
  := \iota_{Y_s}\dot\alpha_s$, where $Y_s$ is the Reeb field of the
  contact form $\alpha_s$.  The solution $X_s$ vanishes on $\partial
  R$, and $X_s$ has a global flow in a small neighborhood of the
  boundary.  Hence one has constructed an $\S^1$-diffeomorphism
  between $R$ and $R^\prime$ that maps the boundary curves onto each
  other, and respects the contact forms close to the boundaries and in
  the neighborhood of the exceptional orbits.
  
  The proof is now finished by applying the Moser trick a second time,
  but now in the interior of the manifold.  The vector field generates
  a global isotopy, because the two contact forms are identical close
  to the boundary components, and the vector field has compact
  support.
\end{proof}

\begin{sphaerenbeispiel}[cont.]
  The Dehn-Euler number $n(R,\gamma)$ is the last invariant that needs
  to be computed to find $(\S^5,\alpha_\pm)$ in the classification
  scheme. The path $\gamma$ can be taken to be $(e^{i\phi},0,0)$ with
  $0\le\phi<2\pi$, and a section in $R =
  \{(z_1,z_2,0)\in\S^5|x_1y_2>x_2y_1\}$ can be found by
  $$
  \sigma:\quad \{z\in\C|\,\mathrm{Im}\,z>0\}\hookrightarrow
  R\subset\S^5,\quad z\mapsto \frac{1}{\sqrt{2+2\abs{z}^2}}
  \big(1+z,z-1,0\big)\;.
  $$
  The boundary of $\sigma$ is composed of two segments
  $1/\sqrt{2}\cdot(e^{i\phi},e^{i\phi},0)$ with $\phi\in[0,\pi]$ and
  $1/\sqrt{2+2x^2}\cdot (x+1,x-1,0)$ with $x\in(-\infty,\infty)$. The
  boundary can be smoothed at the points where the two components
  meet, but this has no effect on the intersection number, because the
  only intersection point of $\partial\sigma$ and $\gamma$ is given by
  $(1,0,0)$, and hence $n(R,\gamma) = \pm 1$. The cross-section $R$
  has opposite orientations for $\alpha_+$ and $\alpha_-$, thus
  $n_+(R,\gamma) = 1$ and $n_-(R,\gamma) = -1$.
  
  The complete set of invariants for $(\S^5,\alpha_\pm)$ is: The
  principal stabilizer is trivial, $\singorb{\S^5}$ has a single
  component that is isomorphic to $\Etwist$, the cross-section is
  $\OpenDisc{2}\times\S^1$, and the Dehn-Euler number $n(R)$
  equals~$\pm 1$.
\end{sphaerenbeispiel}

\begin{brieskornbeispiel}[cont.]
  Above, we already saw that the cross-section of any $W_k^5$ is
  $\S^1$-diffeomorphic to $\OpenDisc{2}\times\S^1$ and
  $\singorb{\left(W_k^5\right)}$ is isomorphic to $\Etriv$ for $k$
  even and $\Etwist$ for $k$ odd.
  
  Now, we will compute $n(R,\gamma)$ for $(W_k^5,\alpha_k)$ and
  $(W_k^5,\alpha_{-k})$.  The curve $\gamma(\phi)$ is given by
  $(e^{i\phi},+i e^{\frac{k}{2}i\phi},0,0)$ with $\phi\in[0,2\pi]$ for
  $k$ even and with $\phi\in[0,4\pi]$ for $k$ odd.
  
  Set $r_0=\abs{z_0}$ and $A =
  \sqrt{2-r_0^2+\sqrt{(2-r_0^2)^2-r_0^{2k}}}$.  The map below is a
  section of~$\overline{R}$
  $$
  \sigma:\quad \Disc{2}\hookrightarrow \overline{R},\quad
  z_0\mapsto \left(z_0, \frac{iz_0^k}{2A}+\frac{iA}{2},
    -\frac{z_0^k}{2A}+\frac{A}{2}, 0\right)\;.
  $$
  The restriction of $\sigma$ to $\partial\overline R$ is
  $\sigma(\phi) = \left(e^{i\phi}, \frac{i}{2}(1+e^{ik\phi}),
    \frac{1}{2}(1-e^{ik\phi}), 0\right)$.
  
  The intersection of $\gamma$ and $\partial\sigma$ is given by the
  equations $2e^{i\frac{k}{2}\phi}= 1+e^{ik\phi}$ and
  $1-e^{ik\phi}=0$, and hence $k\phi=4\pi n$ with $n\in\Z$.  For
  $k=0$, every point of $\partial\sigma$ lies in the curve of marked
  points, but by shifting the section a bit with the $\S^1$-action,
  one obtains $n(R,\gamma)=0$.  For $k$ even, the curve $\gamma$ is
  parametrized by $\phi\in[0,2\pi)$, and there are $k/2$ intersection
  points, for $k$ odd, the curve $\gamma$ closes for
  $\phi\in[0,4\pi)$, and there are $k$ intersection points.
  
  The calculations so far did not depend on the contact form, but one
  can check that $R$ has different orientations for $\alpha_{+k}$ and
  $\alpha_{-k}$.  This changes the orientation of $\partial\sigma$ and
  $\gamma$, but also of $\partial R$, and hence for
  $(W_k^5,\alpha_k)$, we have $n(R,\gamma)=k$, and for
  $(W_k^5,\alpha_{-k})$, we have $n(R,\gamma)=-k$.
  
  The complete set of invariants for $(W_k^5,\alpha_{\pm k})$ is: The
  principal stabilizer is trivial, $\singorb{\bigl(W_k^5\bigr)}$ is
  isomorphic to $\Etwist$ for $k$ odd and to $\Etriv$ for $k$ even,
  the cross-section is $\OpenDisc{2}\times\S^1$, and $n(R)=\pm k$.  In
  particular it follows that the $5$-sphere $(\S^5,\alpha_+)$ in
  Example~\ref{sphaeren beispiel} is equivalent to
  $(W_1^5,\alpha_{+1})$, and $(\S^5,\alpha_-)$ is equivalent to
  $(W_1^5,\alpha_{-1})$.
  
  Note also that every $5$-dimensional simply connected contact
  $\SO(3)$-manifolds with singular orbits is $\SO(3)$-contactomorphic
  to one of the Brieskorn examples $(W_k^5,\alpha_{\pm k})$.  The
  reason is that the orbit space $M/\SO(3)$ of $M$ has to be simply
  connected (\cite{Bredon}), and must have boundary.  Hence $M/\SO(3)$
  is a $2$-disc, and $\singorb{M}$ has a single component. From this
  it follows that the cross-section is isomorphic to
  $\OpenDisc{2}\times\S^1$.  If the principal stabilizer was
  isomorphic to $\Z_2$, then it is easy to show by applying the
  Theorem of Seifert-van Kampen that $\pi_1(M)\cong\Z_2$.  Thus, the
  principal stabilizer has to be trivial, and all cases are covered by
  the $W_k^5$.
\end{brieskornbeispiel}

\subsection{Construction of $5$-manifolds}\label{konstruktion von mannigfaltigkeiten}

In this section, we will construct a manifold $M$ for each of the
possible combination of invariants given in Theorem~\ref{mein
  hauptsatz}.

\subsubsection{$\singorb{M}=\emptyset$}

The classification given in \cite{Kamishima} shows that there is an
$\S^1$-invariant contact structure without Legendrian orbits on any
closed $3$-dimensional contact $\S^1$-manifolds $R$ that does not have
special exceptional orbits or fixed points, and whose (orbifold) Euler
number does not vanish.  The $5$-manifold $M$ is then given by $M\cong
\SO(3)\times_{\S^1}R$, where the circle on $R$ acts with $k$-fold
speed to get the desired stabilizer on $M$.

On the other hand, it follows from Lemma~\ref{singorb ist nullmenge
  von mu} that $0\notin\mu(M)$, and thus $R$ cannot have Legendrian
orbits. It is also clear that $R$ cannot have fixed points.

\subsubsection{$\singorb{M}\ne\emptyset$ and trivial principal stabilizer}\label{dehn twists}

Let $\overline{R}$ be any $3$-dimensional $\S^1$-manifold without
fixed points and without special exceptional orbits, but with
non-empty boundary $\partial R$.  By the requirement that only the
$\S^1$-orbits on the boundary are Legendrian, the contact structure on
$\overline{R}$ is uniquely determined (\cite{Kamishima}).

Over the interior of $\overline{R}$, the $5$-manifold $M^* =
\SO(3)\times_{\S^1}(R-\partial R)$ is a contact $\SO(3)$-manifold.
Now one has to glue in the singular orbits, in such a way as to get
the chosen combination of components of type $\Etriv$ and $\Etwist$
and the Dehn-Euler number $n(R)$.  First we will show how to glue in
the standard model for $\Etriv$; for this, we need to have a standard
form for a neighborhood of $\partial R$.

Let $\sigma$ be any section in $\overline{R}$ that is compatible with
the standard sections around the exceptional orbits. In \cite{Lutz} it
has been shown that any contact form around $\partial R$ is equivalent
to a standard form: Denote the coordinates of a collar $\S^1\times
[0,\epsilon)\times\S^1$ around a boundary component by $(t,r,\phi)$
and let the $\S^1$-action be $e^{i\theta}\cdot(t,r,\phi) =
(t,r,\phi+\theta)$.  Every invariant contact form is up to an
$\S^1$-contactomorphism equal to $dt + r\,d\phi$.  In general the
section $\sigma$ will not be of the form $\sigma(e^{it},r) =
(e^{it},r,1)$ in the collar though, but it is not very difficult to
arrange the model neighborhood in this way. Let $[t]$ and $[\phi]\in
H_1(M,\Z)$ be the classes given by $\S^1\times\{0\}\times\{1\}$ and
$\{1\}\times\{0\}\times\S^1$.  The section $\sigma$ represents an
element $[t]+a[\phi]$, and there is a linear map
$A\in\mathrm{SL}(2,\Z)$ that induces an $\S^1$-diffeomorphism, such
that $\sigma$ represents $[t]$ in the new coordinates.  The contact
form becomes $(1+ar)\,dt + r\,d\phi$, which after dividing by $1+ar$
and rescaling in the $r$-direction can be transformed back into
$dt+r\,d\phi$.  Now by deforming $\sigma$, one obtains a collar for
the boundary where the action, the contact form and the section are
all in standard form.

The standard way of gluing is to consider $\S^1\times T^*\S^2$ with
$\SO(3)$-action on the second factor and with the contact form
$dt+\lcan$.  The cross-section of $\S^1\times T^*\S^2$ looks exactly
like the neighborhood of the boundary components of $\overline{R}$,
which allows us to identify both.  Since the cross-section determines
the $5$-manifold lying over it, this gives a gluing of $\S^1\times
T^*\S^2$ to $M^*$.  In the boundary, the section $\sigma$ and the
curve of marked points are identical, but one can push $\sigma$ a bit
along the $\S^1$-action to avoid having any intersection points.  Thus
the contribution of this gluing to $n(R)$ is zero.

To construct a general $M$, i.e.\ an $M$ with $n(R)\ne 0$ or with
$\Etwist\subset\singorb{M}$, we need to change the construction.

Assume first that we want to glue in a component of type $\Etriv$,
which adds $2c$ to the Dehn-Euler number.  The neighborhood of
$\partial R$ was chosen above to be $\S^1 \times [0,\epsilon) \times
\S^1$ with contact form $dt + r\,d\phi$ and with a section $\sigma$ of
the form $\sigma(e^{it},r) = (e^{it},r,1)$.  The matrix
\begin{align*}
  A &= \left(
    \begin{matrix}
      1 & c \\
      0 & 1
    \end{matrix}
    \right) \in \mathrm{SL}(2,\Z)
\end{align*}
induces a diffeomorphism, which can be isotoped as above to obtain a
new model for the neighborhood of $\partial R$, where $\sigma$
represents the homology class $[t] + c\,[\phi]$, and where the contact
form is still in standard form.  Gluing $\Etriv$ along the
cross-section $\overline{R}$ works again without any problem.  The
intersection number between the section $\sigma$ and the curve of
marked points gives now $c$.

To glue in a component of type $\Etwist$, recall that the
cross-section around $\Etwist$ could be described by $\R\times
[0,\epsilon) \times \S^1 / \!\sim$ with the equivalence relation
$(t,r,e^{i\phi}) \sim (t+1,r,e^{i(\phi + \pi)})$ and contact form
$\alpha = dt + r\,d\phi$.  The curve of marked points was given by
$\{(t,0,1)\}$ and $\{(t,0,-1)\}$.  There is now a diffeomorphism
$\Phi:\, \S^1\times[0,\epsilon)\times\S^1 \to
\R\times[0,\epsilon)\times \S^1/\!\sim,\, (e^{2\pi it},r,e^{i\phi})
\mapsto (t,r,e^{i(\phi + \pi t/2)})$.  The curve of marked points
pulls back to $\{(e^{2\pi it}, 0, e^{-\pi it})\}$, and $\Phi^*\alpha =
(1+ \pi r/2)\,dt + r\,d\phi$, which can be isotoped into standard
form.  The model for the cross-section close to $\Etwist$ and close to
$\partial R$ looks identical, and it is possible to glue both parts.
The Dehn-Euler number $n(R)$ can be arranged in the desired way as
above.

\subsubsection{$\singorb{M}\ne\emptyset$ and principal stabilizer is $\Z_2$}

If the principal stabilizer is isomorphic to $\Z_2$, then all
components of $\singorb{M}$ are equivalent to $\S^1\times\RP{2}$.  The
gluing occurs completely analogous to the way it was done above:
Choose identical charts for a neighborhood of $\partial R$, and for
the cross-section around $\singorb{M}$, and glue along these.

\subsection{Relation between the Dehn-Euler number and generalized Dehn twists}

In this section we want to show that the Dehn-Euler number $n(R)$
counts the number of Dehn twists needed to glue in the singular
orbits.

A \textbf{$k$-fold Dehn twist $\Dehn_k$} on $T^*\S^2$ can be
constructed in the following way.  Write a point in $T^*\S^2$ as
$(\POS,\MOM)\in\R^3\times\R^3$ with $\abs{\POS}=1$ and
$\POS\perp\MOM$.

If one chooses in the map
$$
\Dehn_k(\POS,\MOM) = \Big(\POS\cdot\cos f(\abs{\MOM}) +
\frac{\MOM}{\abs{\MOM}}\cdot\sin f(\abs{\MOM}), \MOM\cdot\cos
f(\abs{\MOM}) - \abs{\MOM}\cdot \POS\cdot\sin f(\abs{\MOM})\Big)
$$
the function $f$ to be $f(r)=r$, then $\Dehn_k$ is just the
standard geodesic flow.  Instead, we will use $f(r)=\pi
k(1+\rho_\epsilon(r))$ with $\rho_\epsilon$ as defined in the proof of
Lemma~\ref{kontaktomorphismus auf R gibt kontaktomorphismus auf M}.
The map is $\SO(3)$-equivariant, and for small $\MOM$, the map is
$\Dehn_k\equiv (-1)^k\id$, while for large $\MOM$, it is
$\Dehn_k\equiv\id$.

The canonical $1$-form transforms like
$$
\Dehn_k^* \lcan = \lcan + \abs{\MOM}\,d\big(f(\abs{\MOM})\big)\;.
$$
This shows that $\Dehn_k$ would be a symplectomorphism of
$(T^*\S^2,d\lcan)$, but not a contactomorphism of $dt+\lcan$ on
$\R\times T^*\S^2$.  It is known (\cite{Seidel}) that $\tau_{2n}$ is
isotopic to $\id$ and $\tau_{2n+1}$ is isotopic to $\tau_1$ (both in
the space of diffeomorphisms with compact support).

The mapping torus
$$
\R\times T^*\S^2/\!\sim, \text{ where $(t;\POS,\MOM)\sim
  (t+1;\Dehn_k(\POS,\MOM))$}
$$
is then diffeomorphic to $V^*\Etriv$ for $k$ even and to
$V^*\Etwist$ for $k$ odd.

The $1$-form
$$
\alpha = dt + \lcan - t\abs{\MOM}\,df
$$
on $\R\times T^*\S^2$ is invariant under the equivalence relation,
and thus projects down onto the mapping torus.  The contact condition
for $\alpha$ gives
$$
0 \ne \alpha\wedge (d\alpha)^2 = dt\wedge (d\lcan)^2 -
2\abs{\MOM}\lcan\wedge dt\wedge df\wedge d\lcan = (1-2\abs{\MOM}^2
f^\prime)\,dt\wedge d\lcan^2\;.
$$
Because $\max f^\prime = c/\epsilon$, by choosing $\epsilon$ small
enough we can assure that $1-2\abs{\MOM}^2f^\prime\ne 0$, and $\alpha$
is then an $\SO(3)$-invariant contact form.

In fact, because the map $\Dehn_k$ used to construct the mapping torus
is the identity far away from the zero-section and because the term
$t\abs{\MOM}\,df$ in $\alpha$ disappears there, it is possible to cut
out a component of $\singorb{M}$ and glue in $\R\times T^*\S^2/\!\sim$
to obtain a new contact $\SO(3)$-manifold.

It only remains to see what effect this has on the integer $n(R)$.

The cross-section $R$ in $\R\times T^*\S^2$ is equal to the one for
the standard contact form
$$
R=\Bigl\{\big(t;(x,y,0),(ry,-rx,0)\big)\Bigr\}/\!\sim\;,
$$
because the last term in $\alpha$ does not change the moment map
($\iota_{X_M}df = \lie{X_M}f$, but $f$ only changes in radial
direction).

To compute the local contribution to $n(R)$, notice that the section
$\sigma(t,r)=(t;(1,0,0),(0,-r,0))$ on $\R\times T^*\S^2$ does not
descend to a continuous section on the mapping torus.  Instead one
could replace $\sigma$ by
$$
\sigma(t,r) = (t;(\cos tf(r), -\sin tf(r),0), (-r\cdot\sin
tf(r),-r\cdot\cos tf(r),0))\;.
$$
Since $\sigma$ remains unchanged far away from the singular orbits,
it extends to the unmodified section, and it is easy to check that
$\sigma$ induces a continuous section on $\R\times T^*\S^2/\!\sim$.

The intersections of $\sigma$ with the curve of marked points is given
by $(\cos tf(0), -\sin tf(0),0)=(\pm 1,0,0)$, i.e. $\cos \pi kt=\pm 1$
and $\sin \pi kt=0$, and then $kt\in\Z$.  There are $k$ points on
$\partial R$, where $\sigma$ intersects the marked set of points.

If $k$ is odd, the boundary corresponds to $\Etwist$.  Then there is
only a single curve of marked points and the contribution of this
boundary to $n(R)$ is $k$.  If $k$ is even, then there are two
disjoint curves of marked points, and there are only $k/2$
intersection points with the first one.  But since for singular orbits
of type $\Etriv$ this number is multiplied by $2$, the contribution to
$n(R)$ is again~$k$.

Thus the Dehn-Euler number $n(R)$ counts the number of Dehn twists
applied at $\singorb{M}$.

All constructions on $\S^1\times\S^2$ in this section are
$\Z_2$-equivariant, and this allows us to build manifolds with
principal stabilizer $\Z_2$ and arbitrary~$n(R)$.

\bibliography{so3-mfld}

\providecommand{\bysame}{\leavevmode\hbox to3em{\hrulefill}\thinspace}
\providecommand{\MR}{\relax\ifhmode\unskip\space\fi MR }
\providecommand{\MRhref}[2]{%
  \href{http://www.ams.org/mathscinet-getitem?mr=#1}{#2}
}
\providecommand{\href}[2]{#2}
\begin{thebibliography}{LMTW98}

\bibitem[Aud91]{Audin}
M.~Audin, \emph{{The topology of torus actions on symplectic manifolds}},
  {Progress in Mathematics, 93. Basel : Birkhäuser Verlag}, 1991.

\bibitem[Bre72]{Bredon}
G.E. Bredon, \emph{Introduction to compact transformation groups}, Academic
  Press, New York, 1972, Pure and Applied Mathematics, Vol. 46.

\bibitem[DK00]{Duistermaat}
J.J. Duistermaat and J.A.C. Kolk, \emph{{Lie groups}}, {Universitext. Berlin:
  Springer}, 2000.

\bibitem[HM68]{Hirzebruch}
F.~Hirzebruch and K.~H. Mayer, \emph{{${\rm O}(n)$}-{M}annigfaltigkeiten,
  exotische {S}ph\"aren und {S}ingularit\"aten}, Lecture Notes in Mathematics,
  No. 57, Springer-Verlag, Berlin, 1968.

\bibitem[Igl91]{Iglesias}
P.~Iglesias, \emph{{Les $SO(3)$-variétés symplectiques et leur classification
  en dimension $4$. (Symplectic $SO(3)$-manifolds and their classification in
  dimension $4$)}}, Bull. Soc. Math. Fr. \textbf{119} (1991), no.~3, 371--396.

\bibitem[KT91]{Kamishima}
Y.~Kamishima and T.~Tsuboi, \emph{{CR-structures on Seifert manifolds}},
  Invent. Math. \textbf{104} (1991), no.~1, 149--163.

\bibitem[LM76]{Lutz_Meckert}
R.~Lutz and C.~Meckert, \emph{Structures de contact sur certaines sph\`eres
  exotiques}, C. R. Acad. Sci. Paris S\'er. A-B \textbf{282} (1976), no.~11,
  Aii, A591--A593.

\bibitem[LMTW98]{Lerman_Meinrenken}
E.~Lerman, E.~Meinrenken, S.~Tolman, and C.~Woodward, \emph{{Non-Abelian
  convexity by symplectic cuts}}, Topology \textbf{37} (1998), no.~2, 245--259.

\bibitem[Lut77]{Lutz}
R.~Lutz, \emph{Structures de contact sur les fibr\'es principaux en cercles de
  dimension trois}, Ann. Inst. Fourier (Grenoble) \textbf{27} (1977), no.~3,
  ix, 1--15.

\bibitem[LW01]{Lerman_Willett}
E.~Lerman and C.~Willett, \emph{{The topological structure of contact and
  symplectic quotients}}, Int. Math. Res. Not. \textbf{2001} (2001), no.~1,
  33--52.

\bibitem[Orl72]{Orlik}
P.~Orlik, \emph{{Seifert manifolds}}, {Lecture Notes in Mathematics. 291.
  Berlin: Springer-Verlag. VIII}, 1972.

\bibitem[Sei98]{Seidel}
P.~Seidel, \emph{{Symplectic automorphisms of $T^*S^2$}}, arXiv
  \textbf{math.DG/9803084} (1998).

\bibitem[vKN]{KoertVan}
Otto van Koert and Klaus Niederkr\"uger, \emph{{Open book decompositions for
  contact structures on Brieskorn manifolds}}, arXiv \textbf{math.SG/0405029}.

\bibitem[Wil02]{Willett}
C.~Willett, \emph{Contact reduction}, Trans. Amer. Math. Soc. \textbf{354}
  (2002), no.~10, 4245--4260 (electronic).

\end{thebibliography}



\end{document}